\documentclass[12pt]{article}

\usepackage{amssymb,amsthm,amsmath, amsfonts,latexsym,epsfig,color, enumerate} 

\usepackage[margin=1in]{geometry}

\newtheorem{theorem}{Theorem}
\newtheorem{lemma}{Lemma}
\newtheorem{proposition}{Proposition}

\newtheorem{example}{Example}

\theoremstyle{definition}

\newtheorem{conjecture}{Conjecture}

\usepackage{hyperref, framed}

\newcommand{\perm}[1]{S_{#1}}

\newcommand{\cho}[2]{{#1 \choose #2}}

\title{Diametric problem for permutations with the Ulam metric (optimal anticodes)}
\author{Pat Devlin\footnote{Swarthmore College, Swarthmore PA, USA \qquad \texttt{pdevlin2@swarthmore.edu}} \and Leo Douhovnikoff\footnote{Swarthmore College, Swarthmore PA, USA \qquad \texttt{ldouhov1@swarthmore.edu}}}

\date{February 29, 2024}

\begin{document}

\maketitle

\renewcommand{\thefootnote}{\fnsymbol{footnote}}
\footnotetext{AMS 2020 subject classification: 05D05, 05A05, 05C35, 05C69}
\footnotetext{Key words and phrases: permutations, isodiametric problem, anticodes, Ulam distance, deletion\\ distance, longest common subsequence, rank modulation, translocations}

\setcounter{footnote}{0}
\renewcommand{\thefootnote}{\arabic{footnote}}

\abstract{
We study the diametric problem (i.e., optimal anticodes) in the space of permutations under the Ulam distance.  That is, let $\perm{n}$ denote the set of permutations on $n$ symbols, and for each $\sigma, \tau \in \perm{n}$, define their Ulam distance as the number of distinct symbols that must be deleted from each until they are equal.   We obtain a near-optimal upper bound on the size of the intersection of two balls in this space, and as a corollary, we prove that a set of diameter at most $k$ has size at most $2^{k + C k^{2/3}} n! / (n-k)!$, compared to the best known construction of size $n!/(n-k)!$.  We also prove that sets of diameter $1$ have at most $n$ elements.
}

\section{Introduction}
One of the most active branches of research within extremal combinatorics and theoretical computer science amounts to the study of finite metric spaces and subsets whose elements satisfy some distance condition.  Most notably, motivated by applications to error-correction, one very popular theme has been to study \textit{codes}---i.e., sets whose elements are all far apart.  For traditional codes, where each point is a binary string, the literature is particularly vast \cite{sloane}, but since the introduction of \textit{rank modulation schemes} \cite{JiangETAL2008} a good amount of recent attention has also been given to codes whose elements are permutations \cite{Barta2017, BuzagloEtzion2015, Han2019PermCodesL1, wang2017}.  In that setting, information is stored as permutations, and various notions of distance are used to model different types of errors that can occur during storage or retrieval.

Dual to the notion of codes are \textit{anticodes}---i.e., sets where any two elements must be close---which in the language of geometry are sets of small diameter.  In addition to having deep connections to coding theory \cite{delsarte}, anticodes are natural objects of study in their own right.  In fact, many celebrated results of extremal combinatorics can be viewed as statements about sets of bounded diameter including the Erd\H{o}s--Ko--Rado theorem \cite{ekr} and---more obviously---Kleitman's theorem \cite{kleitman1966, Ahlswede1998}.  As a more recent example, the question of anticodes of permutations under the Hamming distance, first raised in 1977 by Deza and Frankl \cite{DezaFrankl1977}, was fully resolved in 2010 by an elegant argument of Ellis, Friedgut, and Pilpel \cite{EFP2011}.  See \cite{EllisLifshitz2022, FurediGerbnerVizer2015} for yet more examples as well as discussions of related results.

This paper addresses the metric space of permutations under the \textit{Ulam distance}.  Codes in this space have received a lot of recent attention \cite{FarnoudSkachekMilenkovic2013, hassanzadeh2014multipermutation, farnoud2014, golouglu2015, morePermutationsNotGood2016} and they are related to other combinatorial problems as well \cite{beame2009longest, bukh2014longest, bukh2016twins}.    Our paper addresses anticodes in this space.

\subsubsection*{Ulam distance and our main question}
Let $\perm{n}$ denote the set of permutations on $n$ symbols viewed in one-line notation as strings $\sigma(1) \sigma(2) \ldots  \sigma(n)$. For a permutation $\sigma \in \perm{n}$ and $A \subseteq \{1, 2, \ldots, n\}$, we let $\sigma -A$ denote the string resulting from deleting all the symbols of $A$ from $\sigma$ (and contracting the string to have fewer characters).  For example, we have
\[
45231 - \{1,2\} = 453 \qquad \text{ and } \qquad 15243 - \{5\} = 1243.
\]
In this notation, the \textit{Ulam distance} between strings $\alpha, \beta \in \perm{n}$, can be defined as $d_U (\alpha, \beta) = \min \{ |A| \ : \ \alpha -A = \beta - A\}$.  Said slightly differently, it is the least number of distinct symbols that need to be deleted from $\alpha$ and $\beta$ until they are equal.  This distance, also called the deletion distance, amounts to finding the length of the longest common subsequence, and in fact it satisfies the properties of a metric.  See \cite{Deza2004MetricsOP} for a general discussion of this and other permutation metrics.

The question of determining the expected distance between two uniformly chosen permutations was first raised in 1961 by Ulam \cite{ulam}, and this has since garnered a great deal of attention over the years (see, e.g., the textbook \cite{LISbook} written on this exact subject).  Said geometrically, such results essentially discuss the size of various balls in this metric, with radii approximately $n - 2 \sqrt{n}$ being particularly interesting [as this is the radius required for a ball to occupy a constant proportion of the space].

Instead of studying sets of given radius, we ask how large a set of given \textit{diameter} can be in this space.

\begin{framed}
\noindent \textbf{Main question:} What are the largest possible sets $\mathcal{F} \subseteq \perm{n}$ such that every $\alpha, \beta \in \mathcal{F}$ satisfies $d_U (\alpha, \beta) \leq k$?  Let $f_k (n)$ denote the maximum size of such a set.
\end{framed}

\noindent One particularly natural candidate for a large set with small diameter is:

\begin{example}\label{example:best}
For $A \subseteq \{1,\ldots, n\}$ and $\alpha \in \perm{n}$, consider $\mathcal{F} = \{\gamma \in \perm{n} \ : \ \gamma - A = \alpha - A\}$.  Then $\mathcal{F}$ has size $n! / (n-|A|)!$ and diameter $|A|$.
\end{example}

By the analogy that $\gamma \mapsto \gamma - A$ can be viewed as a type of projection, a set as in Example \ref{example:best} might reasonably be called a cylinder in this space.  This example---as implicitly given in \cite{FarnoudSkachekMilenkovic2013}---shows $\dfrac{n!}{(n-k)!} \leq f_k (n)$, and \cite{ulamSphere2016} computed that when $k$ is even, this set is larger than a ball of radius $k/2$.  This leads us to the conjecture that this construction is in fact best possible.

\begin{conjecture}\label{conjecture:upper bound}
For each $0 \leq k \leq n$, we conjecture $f_k (n) = \dfrac{n!}{(n-k)!}$.  Moreover, if $k < n-2$, the sets attaining this bound must be of the form described in Example \ref{example:best}.
\end{conjecture}

Although the main question above has been recently studied for different metrics on permutations \cite{SchwartzTamo2011, BarLev}, to the authors' knowledge nothing else is known on this question.

\subsection*{Our results}
The first---and most obvious---upper bound on $f_{k}(n)$ is that it must be less than the size of a ball of radius $k$ (it should also be at least the size of a ball of radius $k/2$, but as computed in \cite{ulamSphere2016} this is in fact less than $n!/(n-k)!$).  Since a ball of radius $k$ has size at most $k! {n \choose k}^2$, this gives the bound
\[
\dfrac{n!}{(n-k)!} \leq f_{k}(n) \leq \dfrac{n!}{(n-k)!} {n \choose k}.
\]
When $k$ is small relative to $n$, this implies that $f_{k}(n)$ is between roughly $n^k$ and $n^{2k}/k!$.  Our main result closes this gap considerably---eliminating its dependence on $n$---thereby providing strong asymptotic evidence that Conjecture \ref{conjecture:upper bound} is true.

\begin{theorem}\label{theorem:main-bound}
There is a constant $C \geq 0$ such that for all integers $n$ and $k$ we have
\[
\dfrac{n!}{(n-k)!} \leq f_{k} (n) \leq \dfrac{n!}{(n-k)!} \cdot 2^{k + C k ^{2/3}}.
\]
\end{theorem}

This is obtained as an immediate corollary of the following, which bounds the size of the intersection of two balls in this space and may be of independent interest.

\begin{theorem}\label{theorem:ball-intersection}
There is a constant $C \geq 0$ for which the following is true for all integers $n$ and $k$.  If $\alpha$ and $\beta$ are permutations in $\perm{n}$ with $d_U (\alpha, \beta) = k$, then we have
\[
\Big | \{ \gamma \in \perm{n} \ : \ d_U (\alpha, \gamma) \leq k \textsc{ and } d_U (\beta, \gamma) \leq k \} \Big | \leq 2^{k + C k ^{2/3}} \dfrac{n!}{(n-k)!}.
\]
\end{theorem}

Despite the exponential gap of roughly $2^k$ in Theorem \ref{theorem:main-bound}, it is worth noting that the upper bound via Theorem \ref{theorem:ball-intersection} does not allow for much room for improvement.  Consider for instance $\alpha = 1,2,3, \ldots, n$ and $\beta$ a product of $k$ disjoint transpositions such as
\[
\beta = 2,1,\ 4,3,\ 6,5,\ \ldots,\ 2k, 2k-1,\ 2k+1, 2k+2, 2k+3, \ldots, n.
\]
In this case, $d_U (\alpha, \beta) = k$ and the corresponding intersection of balls as in Theorem \ref{theorem:ball-intersection} has size at least $2^{k - C' k^2 / n} n! / (n-k)!$ (thus when $k = o(n)$, a bound of the form $2^{k + o(k)} n! / (n-k)!$ is best possible in Theorem \ref{theorem:ball-intersection}). 

If desired, Theorem \ref{theorem:ball-intersection} can be generalized to the setting where $d_{U} (\alpha, \beta) = \tau$ is an additional parameter (instead of insisting $\tau = k$), and in fact our proof essentially as written already allows for this.  But we omit this since the details of the argument and result statement become slightly more complicated, and---importantly---doing so did not provide any improved upper bound for Theorem \ref{theorem:main-bound}.  We also note that we made no particular effort to push down the error term $2^{C k^{2/3}}$ since the gap due to the dominant term---namely $2^k$---is unavoidable with this approach.

Turning our attention to other aspects of $f_k (n)$, we first note the following basic properties---included mostly for completeness.
\begin{proposition}\label{proposition:basic properties}
For each $0 \leq k \leq n$, we have
\begin{itemize}
\item[(i)] $f_k (n) \leq f_{k+1} (n)$,
\item[(ii)] $f_k (n) \leq f_k (n+1)$, and
\item[(iii)] $(n+1) f_k (n) \leq f_{k+1} (n+1)$.
\end{itemize}
\end{proposition}
Here, properties (i) and (ii) are immediate from the definition.  Property (iii) follows by a standard sort of `tensoring' trick where given a set $\mathcal{F}\subseteq \perm{n}$ of diameter $k$, we can construct the set
\[
\mathcal{G} = \{ \gamma \in \perm{n+1} \ : \ \gamma - \{n+1\} \in \mathcal{F} \} \subseteq \perm{n+1}.
\]
Then $|\mathcal{G}| = (n+1) |\mathcal{F}|$ and $\mathcal{G}$ has diameter $k+1$.  Property (iii) could be particularly useful in the setting that $n-k = L$ is fixed, in which case we see that $f_{n-L} (n) / n!$ is an increasing (and bounded) sequence in $n$, and Conjecture \ref{conjecture:upper bound} would be equivalent to $\displaystyle \lim_{n \to \infty} f_{n-L} (n) / n! \leq 1/L!$.  Although the authors do not see how to further leverage this observation, we are able to resolve the conjecture in a few special cases:

\begin{proposition}\label{proposition:special cases proven}
For each integer $n$, we have
\begin{enumerate}[(i)]
\item $f_{n} (n) = f_{n-1} (n) = n!$ and $f_0 (n) = 1$,
\item $f_{n-2} (n) = n!/2$,
\item $f_{1} (n) = n$.
\end{enumerate}
\end{proposition}
A complete proof of Proposition \ref{proposition:special cases proven} is given later, but to build some intuition on the problem, let us partially discuss it now.  As before, item (i) follows immediately from the definitions.  Item (ii) is from the simple observation that such a family cannot include both a permutation and its reverse (since those are at distance $n-1$).  Although we were able to come up with several different proofs that $f_1 (n) \leq n$, none seemed particularly well-suited for generalization (even to the case $f_2 (n)$).  In fact, we believe that any human-readable argument resolving $f_{5}(n)$ (for example) would likely generalize considerably.  On the other end, we feel that determining $f_{n-3}(n)$ would also be particularly enlightening.

Prior to discovering Theorem \ref{theorem:main-bound}, we were able to obtain a number of other (strictly worse) results.  Of these, we feel it worth mentioning the following, primarily because its proof is short, and it may perhaps be helpful in generating new approaches.

\begin{proposition}\label{proposition:other bound}
For all $0 \leq k \leq n$ we have,
\[
f_k (n) \leq 4 n (k-1)! \cho{n-1}{k-1}^2 \leq \frac{4 n^{2k-1}}{(k-1)!}.
\]
\end{proposition}

Our main arguments require the following two results, each of which feels too natural (and frankly too easy to prove) to suspect that we are the first to use them.  After a literature search, we were in fact able to find a statement and proof of Lemma \ref{lemma:ulam-graphs} (see e.g., \cite{LCSgraphs}).  We were unable to find a reference for Lemma \ref{lemma:vertex-covers}, but we do not suspect it is new.  Relevant definitions (and a proof of Lemma \ref{lemma:vertex-covers}) are given in Section \ref{section:lemmas}.

\begin{lemma}\label{lemma:ulam-graphs}
For each $\sigma, \gamma \in \perm{n}$, let $G(\sigma, \gamma)$ denote the graph on $\{1, 2, \ldots, n\}$ where $x \sim y$ iff $(\sigma(x) - \sigma(y)) (\gamma(x) - \gamma(y)) < 0$.  Then $\sigma - C = \gamma - C$ iff $C$ is a vertex cover of $G(\sigma, \gamma)$.  In particular, $d_U(\sigma, \gamma)$ is the vertex cover number of $G(\sigma, \gamma)$.
\end{lemma}

The graph $G(\sigma, \tau)$ can also be thought of in terms of \textit{inversion graphs}.  Namely, for each $\gamma \in \perm{n}$, let $I_\gamma$ denote the graph on $\{1, 2, \ldots, n\}$ where $x \sim y$ iff $x < y$ and $\gamma(x) > \gamma(y)$.  Then we have $G(\sigma, \tau) =I_\sigma \oplus I_\tau$, where $\oplus$ denotes the graph symmetric difference [i.e., $e \in I_{\sigma} \oplus I_{\tau}$ iff $e$ is present in exactly one of those inversion graphs].

\begin{lemma}\label{lemma:vertex-covers}
Suppose $G$ is a graph on $V$ with vertex cover number $\tau$.  The number of vertex covers having $m$ vertices is at most $\displaystyle 2^\tau \cho{|V| - \tau}{m-\tau}$.
\end{lemma}

\subsection*{Outline of the paper}
We begin in Section \ref{section:lemmas} where we introduce some common notation and prove Lemma \ref{lemma:vertex-covers}.  We then prove the previously stated Propositions in Section \ref{section:propositions}.  We prove our main results in Section \ref{section:proof-of-ball-bound}, and finally, we conclude in Section \ref{section:conclusion} with a discussion of open problems.

\section{Notation and proof of Lemma \ref{lemma:vertex-covers}}\label{section:lemmas}
Throughout, we focus on the Ulam distance, and we write $d_U$ as $d$.  Recall that a \textit{vertex cover} of a graph is a set $C$ such that every edge of $G$ has at least one endpoint in $C$, and the \textit{vertex cover number} is the size of the smallest vertex cover.  We also use the notation $N(x) = \{y \ : \ y \sim x\}$ to denote the neighborhood of a vertex $x$.

We now turn our attention to a proof of Lemma \ref{lemma:vertex-covers}.

\begin{proof}
Let $S$ be a fixed vertex cover of size $\tau$.  Framing our proof algorithmically, for each vertex cover $C$ of size $m$, we ask the following questions:
\begin{itemize}
\item[(i)] We first ask: what is the set $U := C \cap S$?
\item We then define $\displaystyle T = U \cup \bigcup_{x \in S - C} N(x)$.
\item[(ii)] We then ask: what is $C - T$?
\end{itemize}
Having answered questions (i) and (ii), we claim that this determines $C$ [and in fact that $C = (C-T) \cup T$].  To see this, we need only show that $T \subset C$, which is to say: if $x \in S - C$ then $N(x) \subseteq C$.  In fact, if $y \in N(x)$ but $x \notin C$, then $C$ must meet the edge $x \sim y$ (since $C$ is a vertex cover), so we would need $y \in C$.  Therefore, $T \subseteq C$.

Now we claim that $T$ is a vertex cover.  For this, suppose $u \sim v$ is an edge of $G$.  Since $S$ is a vertex cover, we must have at least one of these vertices in $S$---without loss of generality, say $u \in S$.  If $u \in C$, then $u \in T$ since $u \in C \cap S \subseteq T$.  On the other hand, if $u \notin C$, then $u \in S - C$, so $N(u) \subseteq T$, which means $y \in T$.

Finally, we just need to estimate the number of ways that questions (i) and (ii) can be answered.  Question (i) can be answered in at most $2^{|S|} = 2^{\tau}$ ways.  After this, the set $T$ is determined, and since $T$ is a vertex cover, it has size at least $\tau$.  Thus, question (ii) can be be answered in at most $\displaystyle {{|V|-|T|} \choose {|C| - |T|}} \leq {{|V|-\tau} \choose {m - \tau}}$ ways.
\end{proof}

\section{Proofs of Propositions \ref{proposition:special cases proven} and \ref{proposition:other bound}}\label{section:propositions}
Throughout this section, as a slight abuse of notation, when $\sigma \in \perm{n}$ and $a \in \{1,\ldots, n\}$, we write $\sigma - a$ to mean $\sigma - \{a\}$.

\subsection*{Proof of Proposition \ref{proposition:special cases proven}}
As noted before, the claims in item (i) follow immediately from the definitions, and the fact that $f_{n-2}(n) \leq n!/2$ is because a family of diameter at most $n-2$ cannot have a permutation as well as its reverse (since those are at distance $n-1$).  Thus, we restrict our attention to the claim in item (iii).

\begin{proof}
We need only prove $f_1 (n) \leq n$.  We first outline our argument for clarity's sake.
\begin{itemize}
\item Suppose $\mathcal{F} \subseteq \perm{n}$ is a set of diameter $1$.  For each $1 \leq t \leq n$, let $\mathcal{F}_t \subseteq \mathcal{F}$ denote the set of permutations where the symbol 1 appears in the $t^{\text{th}}$ position.
\item \textbf{Claim 1:} If $|\mathcal{F}_1| \geq 2$, then $\mathcal{F} = \mathcal{F}_1 \cup \mathcal{F}_2$.
\item \textbf{Claim 2:} It is not possible to have $|\mathcal{F}_1| \geq 2$ and $|\mathcal{F}_2| \geq 2$.
\item \textbf{Claim 3:} We have $|\mathcal{F}| \leq n$.
\end{itemize}

\textbf{Claim 1:} Suppose to the contrary that $x,y \in \mathcal{F}_1$ and that we have some $\gamma \in \mathcal{F}_t$ where $t > 2$.  Since $d(x,\gamma) \leq 1$, we know there is some $a$ such that $x - a = \gamma - a$.  But in fact we would need $a = 1$ since otherwise $x - a$ would start with the symbol 1 whereas $\gamma - a$ would not.  Thus we would need $x - 1 = \gamma - 1$ and similarly $y - 1 = \gamma - 1$, but this is impossible since $x - 1 \neq y - 1$ (as $x \neq y$ and both permutations begin with 1).

\textbf{Claim 2:} Now suppose for contradiction $x,y \in \mathcal{F}_1$ and $p,q \in \mathcal{F}_2$.  Say $p= p_1, 1, \ldots$ and $q = q_1, 1, \ldots$.  Since $d(x,p) \leq 1$, there is some $a$ where $x-a = p-a$.  If $a = 1$, then $x = 1, p_1, \ldots$, which would mean that $x- p_1 = p - p_1$.  On the other hand, if $a \neq 1$, then we would need $a = p_1$ since otherwise $x-a$ and $p-a$ would begin with different characters.  Thus, in either case, we must have $x-p_1 = p-p_1$.

Similarly, $x-p_1 = p - p_1 = y-p_1 $ and $x-q_1 = q - q_1 = y-q_1$.

If $p_1 = q_1$, then we would know $p-p_1 = x-p_1 = x-q_1 = q-q_1$, but this would be impossible since $p \neq q$, but in this scenario $p$ and $q$ would both start with the same symbol [the deletion of which would result in equal strings].

Thus, we can assume that $p$ and $q$ have different initial symbols, but we claim this too is impossible.  To see this, if $p - a = q-a$, then we would need $a \in \{p_1, q_1\}$ otherwise they would start with different symbols.  But in that case, one of $p-a$ or $q-a$ would begin with a 1, and the other would have a 1 in its second entry.

\textbf{Claim 3:} Finally, suppose $|\mathcal{F}| > n$.  Then we would know that two permutations in $\mathcal{F}$ would begin with the same symbol---without loss of generality, say this is the symbol 1.  By claims 1 and 2, we would then need that at most one permutation of $\mathcal{F}$ begins with a symbol other than $1$.  And by induction on $n$, we'd have $|\mathcal{F}_1| \leq n-1$, which would mean $|\mathcal{F}| \leq n$ as desired.
\end{proof}

\subsection*{Proof of Proposition \ref{proposition:other bound}}
\begin{proof}
Let $\alpha,\beta \in \mathcal{F}$ be arbitrary and let $\{u, v\}$ be an edge in the graph $G(\alpha, \beta)$ of Lemma \ref{lemma:ulam-graphs} (i.e., say $\alpha$ places $u$ before $v$, but $\beta$ doesn't).  Now suppose $\gamma \in \mathcal{F}$ is arbitrary.  Then we claim that one of the following must hold:
\[
d(\gamma - u, \alpha - u) \leq k-1, \qquad \text{\, \, } \qquad d(\gamma - v, \alpha - v) \leq k-1,
\]
\[
d(\gamma - u, \beta - u) \leq k-1, \qquad \text{or} \qquad d(\gamma - v, \beta - v) \leq k-1.
\]
This is because there must be sets $A, B$ of size $k$ where $\alpha - A = \gamma - A$ and $\beta - B = \gamma - B$.  Thus $\alpha - (A\cup B) = \beta - (A \cup B)$, so $A \cup B$ contains at least one endpoint of the edge $\{u,v\}$.  This means that either $u \in A$ or $v \in A$ or $u \in B$ or $v \in B$ (giving rise to the four cases above).  Moreover, for any $\sigma \in \perm{n}$ and any $w \in \{1,2, \ldots, n\}$ we have
\[
\Big| \{ \gamma \in \perm{n} \ : \ d(\gamma - w, \sigma - w) \leq k-1 \} \Big|  \leq n \cdot |B_{k-1} (n-1)| \leq n \cdot (k-1)! \cho{n-1}{k-1}^2,
\]
where $B_{r} (m) \subseteq \perm{m}$ denotes a ball of radius $r$.
\end{proof}

\section{Proof of Theorem \ref{theorem:ball-intersection}}\label{section:proof-of-ball-bound}
Say $\alpha$ and $\beta$ are in $\perm{n}$ and $d(\alpha, \beta) = k$.  Let $\mathcal{F} = \{\gamma : d(\alpha, \gamma) \leq k \text{ and } d(\beta, \gamma) \leq k\}$. 
We'll bound $|\mathcal{F}|$ by analyzing an algorithm whose outputs include all of $\mathcal{F}$.

\subsubsection*{Process for selecting $\gamma$}
Consider the following process, which---given $\alpha$ and $\beta$---we claim is a way to describe each $\gamma \in \mathcal{F}$ (as well as sets $A$ and $B$ of size $k$ for which $\gamma - A = \alpha -A$ and $\gamma - B = \beta - B$).
\begin{itemize}
\item[(0)] Decide $|A \cup B| = \lambda$.
\item[(1)] Decide $A \cup B$.
\item[(2)] Decide $A$ and $B$.
\item[(3)] For all $x \in A \setminus B$ and $y \in B \setminus A$, decide which comes first in $\gamma$.
\item[(4)] For each $z \in A \cap B$, decide the exact placement of $z$ in $\gamma$.
\end{itemize}
Before further analyzing each of the above steps, let us first prove that going through this process will enable us to uniquely determine $\gamma \in \mathcal{F}$---implying that the number of possible outputs of the above is an upper bound on $|\mathcal{F}|$.

To determine a permutation $\gamma$, note that for all $x, y \in \{1, 2, \ldots, n\}$ we need only decide which symbol appears first.  Having concluded steps (0), (1), and (2), we will have determined both $\gamma - A$ and $\gamma -B$ (since $\gamma - A = \alpha - A$ and $\gamma - B = \beta - B$), which means for each $x ,y \notin A \cap B$, we will know the relative order of $x$ and $y$ unless one symbol is in $A \setminus B$ and the other is in $B \setminus A$.  And thus, after step (3), will will know the relative orders for all $x,y \notin A \cap B$.  Finally, step (4) determines the exact placement of each $z \in A \cap B$ within $\gamma$, which uniquely defines $\gamma$.  This is summarized in the following table.

\begin{table}[h]
\begin{center}
\begin{tabular}{c||c|c|c|c}
& $x \notin A \cup B$ & $x \in A \setminus B$ & $x \in B \setminus A$ & $x \in A \cap B$\\ \hline \hline
$y \notin A \cup B$ & $\gamma-A$, $\gamma-B$ & $\gamma-B$ & $\gamma-A$ & (Step 4)\\ \hline
$y \in A \setminus B$ & $\gamma-B$ & $\gamma-B$ & \textit{(Step 3)} & (Step 4)\\ \hline
$y \in B \setminus A$ & $\gamma-A$ & \textit{(Step 3)} & $\gamma-A$ & (Step 4)\\ \hline
$y \in A \cap B$ & (Step 4) & (Step 4) & (Step 4) & (Step 4)
\end{tabular}
\caption{How the relative order of each $x, y$ is determined}
\end{center}
\end{table}

\subsubsection*{Analyzing this selection process}
To obtain our upper bound on $|\mathcal{F}|$, we bound how many ways each step can be done.

\textbf{Step (0):} First note that $k \leq \lambda \leq 2k$ since $k = |A| \leq \lambda = |A \cup B| \leq |A| + |B| = 2k$.

\textbf{Step (1):} Since $\gamma-A = \alpha - A$ and $\gamma - B = \beta - B$, we know that $\alpha - (A \cup B) = \beta - (A \cup B)$.  Thus, if $G(\alpha, \beta)$ is the graph as in Lemma \ref{lemma:ulam-graphs}, we know $A \cup B$ is a vertex cover of $G(\alpha, \beta)$.  And so by Lemma \ref{lemma:vertex-covers}, there are at most $2^{k} \cho{n-k}{\lambda-k}$ choices in step (1).

\textbf{Step (2):} Having selected $A \cup B$, the number of ways to select subsets $A$ and $B$ each of size $k$ is exactly $\cho{|A \cup B|}{|A|} \cho{|A|}{|A \cap B|} = \cho{\lambda}{k} \cho{k}{2k-\lambda}$.

\textbf{Step (3):} Having already completed steps (0), (1), and (2), we know $A$ and $B$ and we also know $\gamma-A = \alpha - A$ and $\gamma - B = \beta - B$.  For each $x \in B \setminus A$ and $y \in A \setminus B$ we need decide their relative order in $\gamma$.  But since we already know the relative order of all the symbols within $A\setminus B$ and also the relative order within $B \setminus A$, we need only decide how these two ordered lists of symbols interleave.\footnote{This is precisely the same phenomenon as in the well-known ``merge sort" algorithm.}  Thus, step (3) can be done in exactly $\cho{|A \setminus B| + |B \setminus A|}{|A \setminus B|} = \cho{2(\lambda - k)}{\lambda - k}$ ways.

\textbf{Step (4):} This can be done in $\dfrac{n!}{(n-|A \cap B| )!} = \dfrac{n!}{(n-(2k - \lambda))!}$ ways, since this amounts to inserting $|A \cap B|$ symbols into an otherwise sorted list of length $n - |A \cap B|$.

Thus, combining these estimates, we obtain
\begin{eqnarray*}
|\mathcal{F}| &\leq& \sum_{\lambda=k} ^{2k} 2^k \cho{n-k}{\lambda-k} \cdot \cho{\lambda}{k} \cho{k}{2k-\lambda} \cdot \cho{2(\lambda-k)}{\lambda-k} \cdot \dfrac{n!}{(n-(2k-\lambda))!}\\
&=& 2^k \dfrac{n!}{(n-k)!} \sum_{j = 0} ^{k}  {k+j \choose j} {k \choose j} {2j \choose j} \dfrac{(n-k)!}{(n-k +j)!} {{n-k} \choose j}\\
&\leq& 2^k \dfrac{n!}{(n-k)!} \sum_{j = 0} ^{k}  {2k \choose j} {k \choose j} {2j \choose j} \dfrac{1}{j!} \leq 2^k \dfrac{n!}{(n-k)!} \sum_{j = 1} ^{k}  \left( \dfrac{C_0 k^2}{j^3} \right)^{j},
\end{eqnarray*}
for some constant $C_0 > 1$, where the last inequality was obtained using the well-known estimates $\cho{M}{L} \leq (eM/L)^L$ and $M! \geq (M/e)^M$.

Finally, we see that if $F(k,j) = \left( \dfrac{C_0 k^2}{j^3} \right)^{j}$, then $\dfrac{F(k,j+1)}{F(k,j)} = C_0 k^2 \left ( \dfrac{j^j}{(j+1)^{j+1}} \right)^3$ is decreasing, and $F(k,j)$ is maximized when $a \leq j/k^{2/3} \leq b$ for some fixed $0 < a < 1 < b$.  Therefore, we obtain
\[
|\mathcal{F}| \leq 2^k \dfrac{n!}{(n-k)!} \sum_{j = 1} ^{k}  \left( \dfrac{C_0 k^2}{j^3} \right)^{j} \leq 2^k \dfrac{n!}{(n-k)!} \sum_{j = 1} ^{k}  \left(\dfrac{C_0}{a^3} \right)^{b k^{2/3}} \leq 2^{k + C k^{2/3}} \dfrac{n!}{(n-k)!},
\]
as desired. $\qed$

\section{Concluding remarks}\label{section:conclusion}
The most natural direction for future work would be the resolution of Conjecture \ref{conjecture:upper bound}.  This conjecture holds for every value of $f_k (n)$ that we were able to directly compute, although naturally the space of permutations quickly becomes prohibitively large.  A purely combinatorial argument would be particularly pleasing, especially one that relied on techniques such as shifting to directly compare any small diameter set to one of the conjecturally best examples.

As a note of caution, one tempting idea might be to strengthen the main conjecture in the spirit of the Hilton-Milner theorem \cite{hiltonMilner}.  Namely, one might hope that maximal sets which are not of the conjectured type are noticeably smaller.  That said, when $k$ is even, a ball of radius $k/2$ has size roughly $(k/2)! \cho{n}{k/2}^2$, which is annoyingly close to $n!/(n-k)!$, so such an approach---if viable---would need to be done with care.

Another direction for further research might be to determine the best possible bound in Theorem \ref{theorem:ball-intersection}.  Although the example provided following the theorem statement shows that there is limited room for possible improvement, we do not believe that the $2^{C k^{2/3}}$ error term is correct.  In fact, we conjecture that the example discussed yields the maximum possible set size.

Finally, we note a different approach for improving our bound in Theorem \ref{theorem:main-bound} might be to argue that any suitably large family must have two permutations for which the bound given in Theorem \ref{theorem:ball-intersection} is noticeably smaller than the most general case.  For this, we suspect an unavoidable task would be finding permutations $\alpha, \beta \in \mathcal{F}$ for which the graph in Lemma \ref{lemma:ulam-graphs} has appreciably fewer vertex covers than implied by Lemma \ref{lemma:vertex-covers}.

\bibliographystyle{plain}
\bibliography{references.bib}
\end{document}